# INTRANSITIVELY WINNING CHESS PLAYERS' POSITIONS

Alexander Poddiakov
E-mail: apoddiakov@gmail.com

Chess players' positions in intransitive (rock-paper-scissors) relations are considered. Namely, position A of White is preferable (it should be chosen if choice is possible) to position B of Black, position B of Black is preferable to position C of White, position C of White is preferable to position D of Black, but position D of Black is preferable to position A of White. Intransitivity of winningness of chess players' positions is considered to be a consequence of complexity of the chess environment—in contrast with simpler games with transitive positions only. The space of relations between winningness of chess players' positions is non-Euclidean. The Zermelo-von Neumann theorem is complemented by statements about possibility *vs.* impossibility of building pure winning strategies based on the assumption of transitivity of players' positions. Questions about the possibility of intransitive players' positions in other positional games are raised.

## 1. Introduction

Intransitivity of winningness in mathematical objects became widely known beginning from Martin Gardner's mathematical games columns in *Scientific American* [Gardner, 1970, 1974]. This has to do, for example, with non-standard, specially designed dice (intransitive dice) with such figures on their faces that in pair throws die A shows a greater number than die B, die B more often shows a greater number than die C, but die C more often shows a greater number than die A. Accordingly, in order to win, die A has to be chosen in the A-B pair, die B in the B-C pair, and die C in the A-C pair. This accords with the principle of rock-paper-scissors game, and contrasts with the transitive principle of dominance "if A≻B and B≻C then A≻C" seeming universal ("≻" means "is preferable to", "more favorable" etc.).

A note on terminology: the terms "intransitive" and "non-transitive" (e.g., "intransitive dice" and "non-transitive dice") are used for such sets as synonyms in the math literature in spite of some difference between the logical terms "intransitive relation" and "non-transitive relation".

Intransitive sets of playing cards, roulettes, lotteries, etc. may work on the same principle. By now many paradoxical examples of intransitive relation "stochastically greater than" have been invented and many studies in this area have been conducted [Akin, 2019; Bozóki, 2014; Buhler et al., 2018; Conrey et al., 2016; Grime, 2017; Hązła et al., 2020; Hulko & Whitmeyer, 2019; Lebedev, 2019; Polymath, 2017; Trybuła, 1961; Van Deventer, 1992].

Yet what about not stochastic games but deterministic positional games like chess and checkers? Dice are characterized by numbers on their faces, and some dice have such numbers that the dice are intransitive. Chess and checkers are characterized by positions of White and Black. Can these positions be intransitive?

Here the second note on terminology—about term "position"—is necessary. The term has two meanings. On the one hand, it means one player's position (e.g., in comparison with the other player's position). Some examples of use of the term in this meaning are the following.

David Brodsky [2018] entitles his post *Winning Equal Positions: An Imperfect Job* and writes: "Both players' positions are fairly solid. White's rook on h1 isn't doing anything for the moment, and his superfluous knights aren't awe-inspiring, but he can improve his pieces". "You need to learn to find weaknesses in your opponent's positions" [Van Apeldoorn, 2014]. "If the idea of balance was sought of as an old fashion scale, like the scales of justice, both players' position would hang equally in relation to one another. However, if one player has better development the scale will tilt in his or her favor" [Patterson, 2014].

On the other hand, "position" means a whole including both players' positions (positions of White and Black): "The position on the board looks very interesting" (it is not about one player's position but about interaction of both players' positions). Here one should cite the following statements concerning transitivity of positions in the second meaning (as wholes).

"Current schemes for machine evaluation of a chess position are based on two assumptions: First, for any two positions, it is possible to decide which is more favorable and second, that this relation is transitive. In consequence of these assumptions, a number that serves as a measure of worth can be assigned to each possible chess position" [Atkinson, 1998, p. 38].

"The construction of a checkers endgame database is simply the computation of a transitive closure. Each position is a member of either the set of wins, losses or draws. Once computed, the classification of a database entry represents perfect knowledge as to the theoretical value of that position" [Lake et al., 1993, p. 3].

Account of perfect values based on playout values has been used to build a common evaluation function for playing Chinese Checkers with two or more players [Sturtevant, 2015].

To distinguish between these meanings of term "position", we will use terms "players' positions" (positions of players, positions of sides) *vs.* "whole positions".

A clarification is in order which is important for practical play but not for theory. Two above-cited works consider sets only of such pairs of whole positions transition between which according to the rules of the game is possible during playout. Based on this approach, checkers has been solved [Schaeffer et al., 2007]. Without questioning the postulated transitivity of whole positions and results achieved ("checkers is solved") in a set of positions, arising of which is possible in playout games, we can tackle the question of whether it is possible to calculate perfect theoretical values of each player's position in the whole set of all positions in the first meaning.

The main question of the article: are transitive ordering of all positions of players and account of theoretical, perfect values of players' positions possible in chess and checkers?

The answer is: it is impossible.

## 2. Examples of Intransitively Winning Chess and Checkers Players' Positions

Let us consider specially constructed chess positions chains of which cannot arise in a game unfolding according to classical rules, but which are not forbidden by these rules.

Let us consider the following four positions: position A of White, position B of Black, position C of White, and position D of Black (Fig. 1) [Poddiakov, 2016, p. 48].

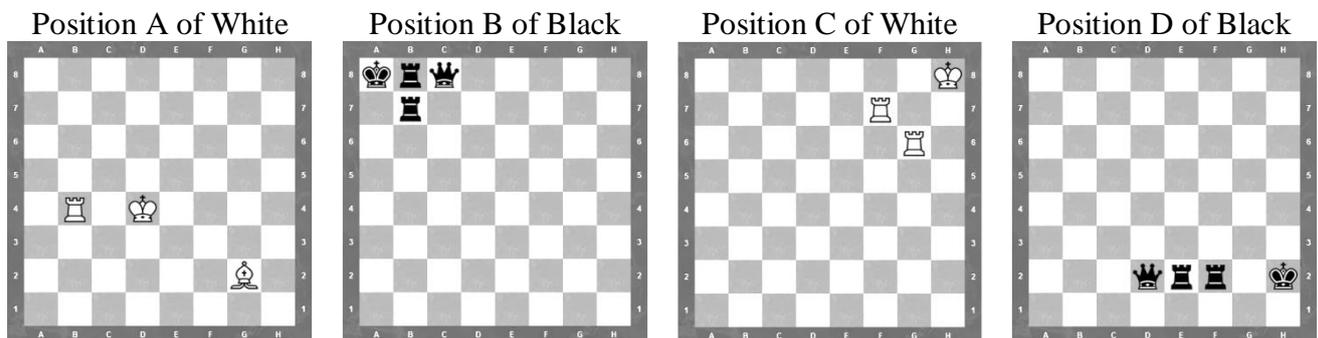

Fig. 1. Two chess positions of White and two positions of Black to be superposed for demonstration of intransitivity of their winningness [Poddiakov, 2016].

Now let us check their pairwise superpositions. Let us lay positions
- A and B
- B and C
- C and D
- A and D
      on a chess board (Fig. 2).

White start in all the compositions in accordance with the rules of chess problem composing.

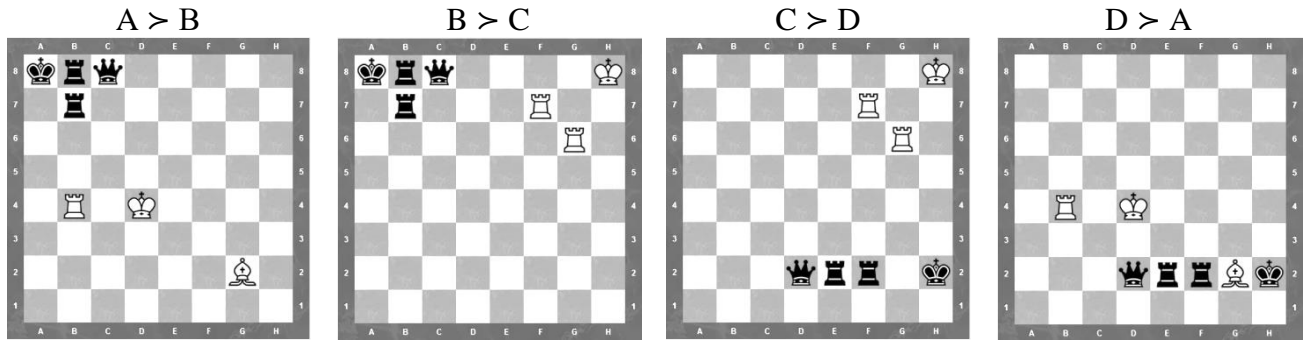

Fig. 2. The chess players' positions superposed for demonstration of intransitivity of their winningness [Ibid.].

One can see that:
- position A of White is preferable (it should be chosen if choice is possible) to position B of Black;
- position B of Black is preferable to position C of White;
- position C of White is preferable to position D of Black;
- but position D of Black is preferable to position A of White
(like in rock-paper-scissors game).

What is the minimum quantity of pieces providing with opportunity of intransitively winning players' positions? Alexander Filatov [2017] has designed a minimalist and symmetrical intransitive chess players' positions (a lesser number of pieces is impossible) (Fig. 3).

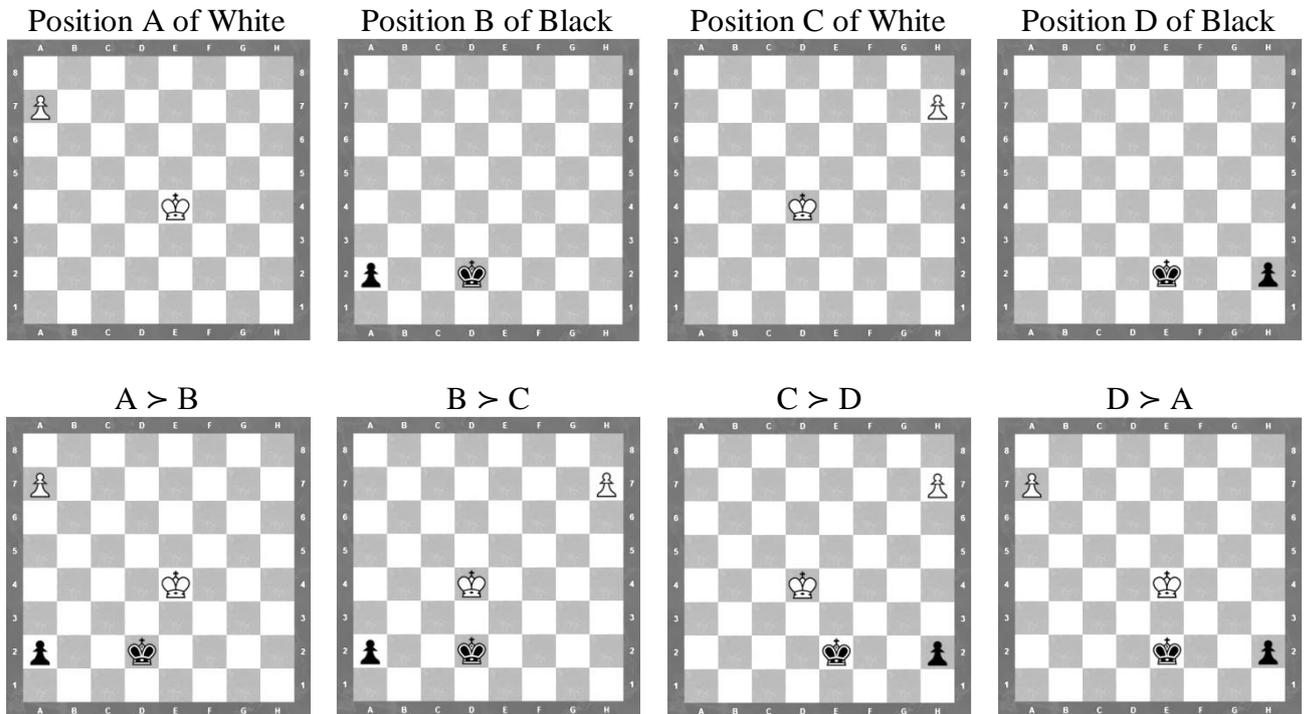

Fig. 3. Filatov's minimalist and symmetrical structure of intransitively winning chess players' positions [2019].

Various chess problems can be designed based on such intransitivity. Grigory Popov, a chess composer, has designed the following problem (Fig. 4). For a possible solution see [Popov, 2021].

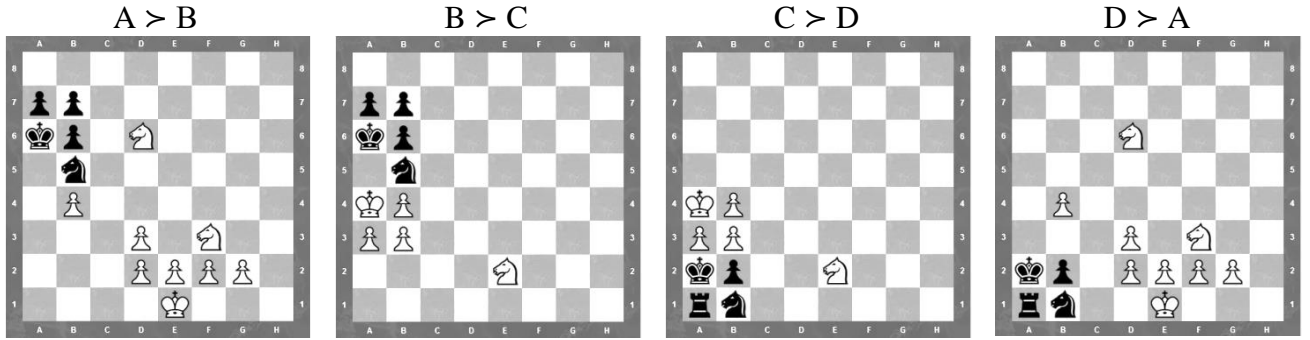

Fig. 4. Popov's chess problem based on intransitivity [Popov, 2021].

## 3. How Frequent are Intransitive Cycles of Chess Players' Positions?

Alexander Filatov [2017], using my example presented above, has shown that the number of intransitive positions in chess is astronomically large, and that intransitive chess chains can be of astronomical length (so it is futile to try to build the longest one).

What is the share of intransitive positions in the total body of positions?

No answer has been given, but we propose an approach to solving the question.

Since there exist chess endgame databases, which make it possible to calculate the outcome of any endgame for 3-7 pieces (review and comparison of the databases are given by Peterson [2018]), we can do the following.

1. Consider 4-7-piece endgames varying the length of position chains (4, 6, 8 positions). As it is difficult to go over all of them, let us fall back on the Monte-Carlo method to generate chains: we randomly generate position A of White (two pieces—the King and some other piece), position B of Black of two pieces, position C of White of two pieces and position D of Black of two pieces (for a chain of length 4). We see if an intransitive cycle of winningness is formed in this chain. Repeat the same for 4 new generated positions. Do it multiple times. Count the number of intransitive cases among all the cases considered. A similar procedure has been used to evaluate the share of intransitive chains for N-sided dice [Conrey et al., 2016]).

Counterintuitively, but by analogy with the results of (Ibid.), with a large number of positions (A, B, C, D, E, F; A, B, C, D, E, F, G, H; etc.) the probability of encountering an intransitive chain does not radically diminish, remaining at least a value of the same order as the probability of drawing a transitive sequence and may even larger.

2. Compare the results for endgames with varying numbers of pieces and unvarying length of chains. Again, counterintuitively, the probability of "drawing" an intransitive chain may grow with the growth of the number of pieces in chains of equal length.

Intransitive checkers players' positions are also possible. See an example below (Fig. 5).

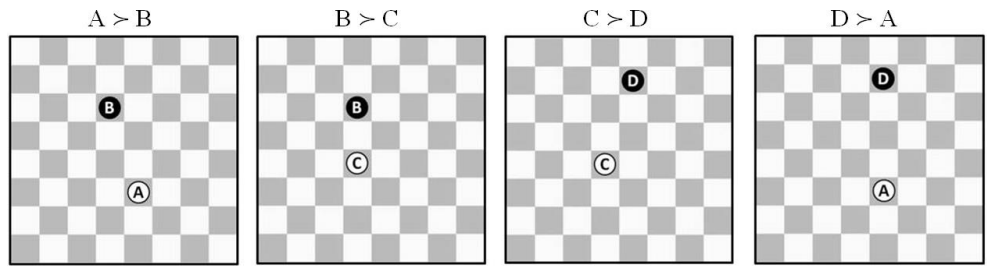

Fig. 5. Zhurakhovsky' intransitive checkers players' positions, in context of Poddiakov's intransitive chess positions: A ≻ B, B ≻ C, C ≻ D, D ≻ A.

Using the methods analogous to those used by Filatov, it is apparently possible to show that the number of intransitive chains in checkers is also very large.

The intransitivity of players' positions shown above hardly matters for practical play. Let us repeat that the intransitively winning players' positions exist in the space of all positions (and not in an ordinary played game) and are only revealed when superposed on one board artificially. But intransitivity of players' positions is important for theoretical assessment of (im)possibility of transitive ordering of all positions and for calculating the perfect value of each of them.

## 4. Intransitivity of Chess Players' Positions Is A Consequence of Complexity of Chess Environment

Intransitivity of chess players' positions is a previously unknown property of chess. No one initially sought to create conditions for the building of intransitive chess players' positions (including chains of astronomical length). This non-evident property was an unplanned result of a complicated environment created by humans.

I think that intransitive players' positions become impossible in simpler environments—e.g., on small chess boards (3×3? 4×4?). If so, it is interesting to determine the minimum size of the board on which intransitivity gets already possible.

And what about shapes of boards? Are transitive chess players' positions possible on toroidal and cylinric boards?

Oleg Yarygin's hypothesis (informal communication): there exists the number of moves between chess players' positions that must be exceeded before intransitivity becomes possible. It would be interesting to determine the limit beyond which intransitivity becomes possible and find out whether it can be a measure of the game's complexity.

Intransitivity of players' positions may exist in other complex positional games, for example, in Chinese Checkers, Go, Reversi etc. Perhaps, the number and ratio of intransitive players' positions among all positions, and their diversity can serve as measures of complexity of a positional game.

Having shown that intransitive players' positions are possible in complicated games, let us now look at a simpler game with a *transitive order* of players' positions' winningness and the possibility of calculating their perfect values. We would need this for further comparison.

## 5. The Magicians: A Simpler Game with Transitive Positions of Pieces Only

I have designed The Magicians to study children's thinking and problem solving. In some aspects, The Magicians are similar to Reversi (see info about this game in [MacGuire, 2011]) but is much simpler than it.

The material of the game are similar cards. Each card has the picture of a good magician (a smiling rosy face) on one side and a bad magician (an angry blue face) on the other side (like two-faced Janus).

The cards are arranged in two horizontal rows on a marked board with two rows and some number of cards in cells of each row (this number can vary in different versions of the game).

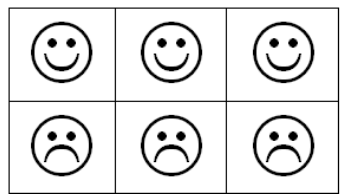

Fig. 6. Composition of the cards in Problem "Three good magicians against 3 bad ones".

*Rules of the game*
A player can play with oneself, and this version of the game is similar with solitaire. Game of two players is possible too, but it is not very interesting for adults because of simplicity.

By shifting the cards, bad magicians can be turned into good ones and vice versa. The aim of the game is to turn all bad ones into good ones. The move consists in swapping one card from the upper row for one card in the lower row. A swap of any card from the lower row and any card in the upper row is possible independently of their positions in the rows (swaps can be vertical and across from one another, it does not matter in this game). Swapping within the row is forbidden. Moves of cards to empty places outside the marked board are forbidden

If, as a result of the swap a magician swapped finds oneself between two wrong ones (i.e. a bad one between two good ones or a good one between two bad ones) s/he too becomes a "wrong one") (the card is turned with the other side up). "You've got surrounded and the same as surroundings" (a mnemonic rule for children).

In all other variants of neighborhood no transformation occurs. It has to be stressed that only magicians who are new to a position are liable to be transformed. The rest do not change even if they get hostile neighbors as a result of the swap near them.

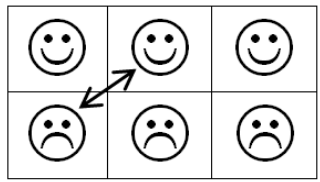

Fig. 7. A useful move (swap) in Problem "Three good magicians against 3 bad ones". The bad magician swapped from the lower row is turned into a good one (because s/he gets between 2 good ones in the upper row), but the good magician swapped from the upper row does not change (because s/he does not get 2 "enemies"-neighbors in the lower row). So, we will have 4 good magicians and 2 bad ones—instead of 3 against 3.

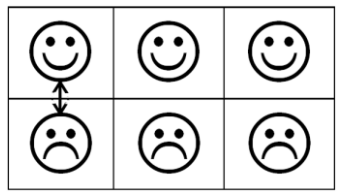

Fig. 8. A possible but worsening move (swap) in Problem "Three good magicians against 3 bad ones". Both magicians do not change (because they do not get 2 "enemies"-neighbors after the swap, but working pairs disappear, and it is bad (a working pair is two friendly magicians separated by one cell in which they convert foes under favorable circumstances, like on a conveyor belt).

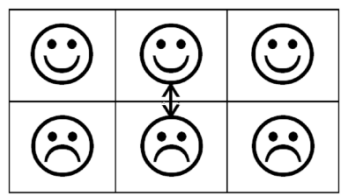

Fig. 9. Another possible, not worsening but useless move (swap) in Problem "Three good magicians against 3 bad ones". After the swap, both magicians are turned into their "enemies" by new hostile neighbors, and the situation does not change.

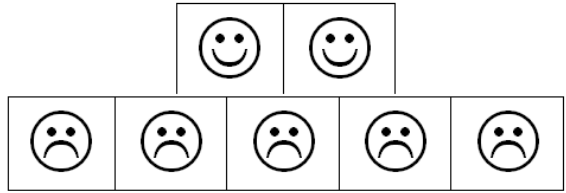

Fig. 10. Problem "Two good magicians against 5 bad ones". Preliminary swaps without the magicians' transformations are necessary here, and it is difficult for most children.

I have devised a computer version of the game for which I have written an algorithm involving an active virtual enemy (the chief bad magician who plays for the bad ones).

This algorithm takes into account only: the position of friendly magicians in a row (the closer the magician is to the middle of the row or to creation of a working pair the higher his value) and already existing working pairs.

It is a simple game: the assessment of the position of "friends" does not include the position of the "strangers", *all the positions (down to symmetrical ones) have their unique values, and there is no intransitivity of players' positions*.

An example of the transitive order of a player' positions in one row for board $3 \times 2$ (like on Fig. 7) is the following.

OOO ≺ OXO ≺ XOO ≺ XXO ≺ XOX ≺ XXX,
where X means a good magician, and O means a bad one.

Here all the combinations in one row (without symmetrical ones) are listed.

XOX ≻ XXO, because XOX contains a "working pair" of good magicians able to turn bad magicians—in contrast with XXO.

XOO ≻ OXO, because XOO needs only one move (swap) to be turned into XOX in which good magicians form the "working pair". OXO needs at least 2 moves for it.

A sum of values in both rows is counted, and the swap leading the maximum sum is chosen. In case of equal sums for symmetrical positions (the algorithm does not identify symmetrical positions as symmetrical), a random choice of the appropriate swaps is made. Transitivity works here too.

Thus, there exist two types of deterministic positional zero-sum games.

1. Games without intransitively winning players' positions. All the positions of players can be ordered transitively, and perfect values can be built for them (like in The Magicians).

2. Games with intransitively winning players' positions (chess, checkers). Positions of each side in such games do not lend themselves to absolute evaluation and absolute rating (there can be no perfect values for all the positions).

## 6. Conclusion

Intransitive objects are of various types and include not only stochastic ones but also deterministic objects—players' positions in abstract strategy games.

Intransitivity of winningness of players' positions is a consequence, a by-product of complexity of the games. The complexity of the chess (checkers) environment makes intransitive players' positions possible—in contrast to simpler games.

Peter C. Fishburn, the author of a theory of decision making without axiom of transitivity, wrote that rejection of intransitivity is analogous with "rejection of non-Euclidean geometry in physics", which "would have kept the familiar and simpler Newtonian mechanics in place, but that was not to be" [Fishburn, 1991, p. 117].

In turn, we have shown that, for chess and checkers, any Euclidean metric of players' positions across their whole set is impossible because some positions form intransitive cycles. The space of mutual relations between winningness of chess players' positions is non-Euclidean. There are no perfect values (numbers in any absolute rating) of all the positions of players in Euclidean space.

Let us consider it in more detail.

Perfect values of players' positions in case of their intransitive cycling should be equal to one another (each position beats another position and is beaten by the third one—like mutual beating of intransitive dice which do not have any absolute winners and losers). The Euclidian distance between them should be equal to zero (there should be no any preferences between them).

Yet if position A beats position B then the value of A should be higher than the value of B (the distance between the positions does exist), and so on, and if the last member of the sequence beats the first one then the value of this last member should be higher than the value of the first one (the distance between the positions does exist).

Thus, perfect values (numbers) should describe intransitive loop "$N_A>N_B$, $N_B>N_C$, $N_C>N_D$, $N_D>N_A$", where $N_A$, $N_B$, $N_C$, and $N_D$ – perfect values for positions A, B, C, and D, respectively. Such perfect values are impossible.

Thus, the existence of intransitive players' positions means that a number that serves as a measure of worth (preference, winningness) cannot be assigned to each possible chess (checkers) position of a player. There are no perfect values of all the positions of players in Euclidean space.

The number of intransitive players' positions in chess is huge, but one example, like in checkers, suffices to prove that perfect values of positions are impossible.

What is called the Zermelo-von Neumann theorem consists in that all positional chess-like games with perfect information have solutions in pure strategies [Csákány, 2002; Kolchin].

For my part, the following should be added.

For some of these positional games a solution is possible on the basis of the assumption on the transitivity of the winningness of players' positions in the whole set of the positions (like in Magicians). But for some games such a blanket solution is impossible (because they have intransitively winning players' positions).

For both types of games, intransitivity of winningness of positions as wholes is not questioned. We have discussed intransitivity of winningness of players' positions only and some of its consequences listed above.

## Acknowledgments

I would like to thank David Silver, Yaroslav Shitov, Oleg Yarygin, and an anonymous reviewer for discussion of the manuscript and important comments and recommendations.